\documentclass[12pt]{article}
\usepackage{amsmath}
\usepackage{amssymb}
\usepackage{enumitem}

\newtheorem{theorem}{Theorem}

\newtheorem{claim}{Claim}
\def \bthm {\begin{theorem}}
\def \ethm {\end{theorem}}
\newtheorem{defi}{Definition}

\newtheorem{prop}{Proposition}
\newtheorem{rmk}{Remark}
\newtheorem{cor}{Corollary}
\newtheorem{exa}{Example}
\newtheorem{conj}{Conjecture}
\newtheorem{prm}{Problem}

\def \nin {\noindent}
\def \bsk {\bigskip}
\def \msk {\medskip}

\def \pf {\nin{\bf Proof. }}
\def \qed {\hfill $\Box$}

\def \cir {\chi_1}
\def \omi {\omega_1}
\def \ali {\alpha_1}
\def \thi {\theta_1}
\def \smin {\setminus}
\def \es {\varnothing}

\def \cG {\mathcal{G}}
\def \cP {\mathcal{P}}
\def \cT {\mathcal{T}}
\def \tw {\textrm{tw}}
\def \tb {\mathsf{Tab}}

\def \ve {\varepsilon}
\def \vp {\varphi}
\def \np {{\sf NP}}
\def \npc {\np-complete}
\def \nph {\np-hard}

\begin{document}

\title{The robust chromatic number of graphs}
\author{G\'abor Bacs\'o$^a$\and
    Bal\'azs Patk\'os$^b$\and
  Zsolt Tuza$^{b,c}$  \and
  M\'at\'e Vizer$^d$ \\
  \small $^a$ Institute for
  Computer Science and Control \\
\small $^b$ Alfr\'ed R\'enyi Institute of Mathematics \\
\small $^c$ University of Pannonia \\
\small $^d$ Budapest University of Technology and Economics
  }
\date{
}
\maketitle

\sf

\begin{abstract}
A 1-removed subgraph $G_f$ of a graph $G=(V,E)$ is obtained by
 
\begin{itemize}
  \item[$(i)$] selecting at most one edge $f(v)$ for each vertex $v\in V$, such that $v\in f(v)\in E$
(the mapping $f:V\to E \cup \{\varnothing\}$ is allowed to be non-injective), and
  \item[$(ii)$] deleting all the selected edges $f(v)$ from the edge set $E$ of $G$.
 \end{itemize}

Proper vertex colorings of 1-removed subgraphs proved to be a useful tool for earlier research on some Turán-type problems.

In this paper, we introduce a systematic investigation of the graph invariant 1-robust chromatic number, denoted as $\cir(G)$. This invariant is defined as the minimum chromatic number $\chi(G_f)$ among all 1-removed subgraphs $G_f$ of $G$. We also examine other standard graph invariants in a similar manner.
\end{abstract}

\section{Introduction}

We consider finite simple graphs $G=(V,E)$, without loops and
 multiple edges. 
Sometimes the notation $V(G)$ and $E(G)$ will also be used.
For the chromatic number, clique number, independence number,
 and clique covering number
 we use the standard notation $\chi(G)$, $\omega(G)$, $\alpha(G)$,
 and $\theta(G)$, respectively.

\begin{defi}
$(i)$
A\/ \emph{1-selection} in a graph\/ $G=(V(G),E(G))$ is a mapping\/
 $f:V(G)\to E(G) \cup \{\varnothing\}$ such that\/ $v\in f(v)$ holds for all\/ $v\in V(G)$ with $f(v)\neq \varnothing$.
The graph\/ $G_f$ with vertex set\/ $V(G_f)=V(G)$ and edge set
 $$
   E(G_f) := E(G) \setminus f(V(G)) 
 $$
  is termed a\/ \emph{1-removed subgraph} of\/ $G$.

$(ii)$
A graph is said to be \emph{quasi-unicyclic} if
 each of its components is a tree or a unicyclic graph.
\end{defi}

A transparent representation of a 1-selection $f$ can be given
 by a directed graph $D(G,f)$ whose vertex set is $V(G)$, and
 for each $v\in V(G)$ with a non-empty image the selection $f(v)=vw$ is represented
 as the arc $(v,w)$, which is oriented from $v$ to $w$.
Hence directed cycles of length $2$ may also occur.
According to the definitions, all vertices have out-degree at most $1$
 in $D(G,f)$.
For this reason the underlying undirected graphs of $D(G,f)$ for 1-selections $f$
 are quasi-unicyclic.

We introduce the following graph invariants concerning
 1-removed subgraphs.

\begin{defi}   \label{d:param}
For a graph\/ $G$,
 \begin{itemize}
  \item the\/ \emph{$1$-robust chromatic number} of\/ $G$ is
 $
   \cir(G) := \min_f \chi(G_f),
 $
  \item the\/ \emph{$1$-robust clique number} of\/ $G$ is
 $
   \omi(G) := \min_f \omega(G_f),
 $
  \item the\/ \emph{$1$-robust independence number} of\/ $G$ is
 $
   \ali(G) := \max_f \alpha(G_f),
 $
  \item the\/ \emph{$1$-robust clique covering number} of\/ $G$ is
 $
   \thi(G) := \max_f \theta(G_f),
 $ 
  \item the\/ \emph{$1$-robust chromatic index} of\/ $G$ is
 $
   \cir'(G) := \min_f \chi'(G_f),
 $
 \end{itemize}
 where\/ $\min$ and\/ $\max$ are taken over all 1-selections\/ $f$
 on\/ $G$.
\end{defi}

The following proposition collects some basic properties of the
 1-robust chromatic number.
The proofs are immediate from the definitions.

\begin{prop}   \label{p:basic}
$(i)$ The value\/ $\cir(G)$ of a graph\/ $G=(V,E)$ is equal
 to the minimum number\/ $k$ of vertex classes in a partition\/
 $V=V_1\cup\cdots\cup V_k$ such that each\/ $V_i$ induces a
 quasi-unicyclic subgraph in\/ $G$.

$(ii)$ A graph\/ $G$ satisfies\/ $\cir(G)=1$
 if and only if it is quasi-unicyclic.
In particular, every tree has\/ $\cir=1$.

$(iii)$ If\/ $G=(V,E)$ does not have any tree components, then
 in computing\/ $\cir(G)$ one may restrict attention to
 1-selections\/ $f$ that are injective, i.e.\/ $f(v)\neq f(v')$
 for any two distinct\/ $v,v'\in V$, without loss of generality.
\end{prop}

Due to $(i)$, $\cir$ corresponds to a weakening of the condition
 that defines ``vertex arboricity'', as in the latter only a subfamily of
 1-selections is allowed; cf.\ Section~\ref{ss:defs}
  and later Proposition \ref{p:arb}.

In fact, the above notions can be put in a more general setting.

\begin{defi}
For a non-negative integer $s$, an \emph{$s$-selection} on a graph $G=(V, E)$ is defined as a function $f: V \to 2^E$ such that $f(v)\subseteq E(v)$ and $|f(v)|\leq s$  where $E(v)$ refers to the set of edges incident with vertex $v$. Using this definition, one can introduce various graph parameters like $\chi_s(G)$, $\omega_s(G)$, $\alpha_s(G)$, etc. as defined in Definition~\ref{d:param} by taking the minimum or maximum value over all $s$-selections.
\end{defi}

With this formalism the chromatic number, the clique number,
 the independence number, and the clique covering number of $G$
 may be viewed as $\chi(G)=\chi_0(G)$, $\omega(G)=\omega_0(G)$,
 $\alpha(G)=\alpha_0(G)$, and $\theta(G)=\theta_0(G)$, respectively.
Then standard inequalities generalize as follows.

\begin{prop}
For every graph\/ $G$ and every integer\/ $s\geq 0$ we have
 $$
   \chi_s(G) \geq \omega_s(G) \qquad \mbox{\rm and} \qquad
     \chi_s(G) \geq \frac{|V(G)|}{\alpha_s(G)} ,
 $$
 moreover
 $$
   \theta_s(G) \geq \alpha_s(G) \qquad \mbox{\rm and} \qquad
     \theta_s(G) \geq \frac{|V(G)|}{\omega_s(G)} .
 $$
\end{prop}

\paragraph{Simplified terminology.}

In the sequel we concentrate on the case of $s=1$, leaving the
 larger values of $s$ for later research.
For this reason, we will just write ``robust'' instead of
 ``1-robust'' for each of the parameters
  $\cir$, $\omi$, $\ali$, $\thi$, $\cir'$.

\begin{exa}   \label{ex:kn}
If\/ $G$ is the complete graph\/ $K_n$, which has\/ $\chi(G)
 = \omega(G) = n$, then\/ $\cir(G)=\omi(G)=\lceil n/3 \rceil$.
The lower bound\/ $\omi(G)\geq\lceil n/3 \rceil$ is a direct consequence of Tur\'an's theorem, as
 at most\/ $n$ edges are removed using a 1-selection.
The upper bound\/ $\cir(G)\leq\lceil n/3 \rceil$ is easily seen by
 splitting the vertex set into\/ $\lceil n/3 \rceil$ disjoint sets\/
  $V_1,\dots,V_{\lceil n/3 \rceil}$ of sizes at most\/ $3$,
 and removing all edges inside each\/ $V_i$.
\end{exa}

\begin{exa}   \label{ex:kttt}
Let\/ $t>k\geq 2$ be any integers.
If\/ $G$ is the complete\/ $k$-partite graph\/ $K_{t,\dots,t}$,
 then of course\/ $\chi(G) = \omega(G) = k$.
But we also have\/ $\cir(G)=\omi(G)=k$.
Indeed, $G$ contains\/ $t^k$ copies of\/ $K_k$, and each edge is
 contained in exactly\/ $t^{k-2}$ copies of\/ $K_k$.
The number of vertices is\/ $kt$, hence by removing that many edges,
 no more than\/ $kt \cdot t^{k-2}=kt^{k-1}$ copies of\/ $K_k$ can be destroyed
 and the clique number\/ remains\/~$k$.
\end{exa}

Below we shall see that the conclusion $\cir(G)=\omi(G)=k$ is
 valid also for $t=k$ if $k\geq 3$.
(This is not the case if $k=2$ because $\cir(K_{2,2}) = \cir(C_4) = 1$.)

\subsection{Motivation and earlier results}

The robust chromatic number $\cir$ was introduced in \cite{PTV}
 as a useful tool to derive estimates on a Tur\'an-type
 extremal problem on graphs with edges assigned sets of integer vectors. The paper established the following results on $\cir$ for complete multipartite graphs and random multipartite graphs.

\bthm   \label{t:PTV}
$(i)$ \cite[Proposition 2.6.]{PTV}
The complete tripartite graph\/ $K_{r,s,t}$ with\/ $1 \le r \le s \le t$
 and\/ $t\ge 2$ satisfies\/  $\cir(K_{r,s,t}) = 2$ if and only if\/  $r \le 2$;
 otherwise\/ $\cir(K_{r,s,t}) = \chi(K_{r,s,t}) = 3$.
 (If\/ $t=1$, for\/ $K_{1,1,1}$ we have\/ $\cir(K_{1,1,1}) = \cir(C_3) = 1$.)

$(ii)$ \cite[Theorem 1.10.]{PTV}
Let\/ $K(m, r, p)$ denote the probability space of all labeled\/
 $r$-partite graphs with each partite set having size\/ $m$, where
  any two vertices in different parts are joined with
   probability\/ $p$, independently of any other pairs.
If\/ $p = \omega (m^{-1/\binom{r}{2}})$, then\/
 $\cir(K(m, r, p)) = r = \chi(K(m, r, p))$ with probability
 tending to\/ $1$ as\/ $m$ tends to infinity.

$(iii)$ \cite{PTV}
A bipartite graph\/ $F$ has\/ $\cir(F)=2$ (i.e.,\/
 $\cir(F)=\chi(F)$) if and only if it contains a component with
 more edges than vertices.
\ethm

\subsection{Standard definitions and notation}   \label{ss:defs}

Beside $\alpha,\omega,\chi,\theta,\chi'$ which already
 occurred above, we use the standard notation $\delta(G)$ for
 minimum vertex degree and $\Delta(G)$ for maximum vertex degree.
Also, for two graphs $G$ and $H$ we write $G \oplus H$ to denote
 their complete join if $G$ and $H$ are vertex-disjoint; and
  $G\cup H$ will denote the graph with vertex set $V(G)\cup V(H)$
 and edge set $E(G)\cup E(H)$, where $V(G)\cap V(H)=\es$
 may or may not hold.
The \emph{lexicographic product} of $G$ and $H$, denoted by
 $G\circ H$, has vertex set $V(G\circ H) = V(G) \times V(H)$;
 two vertices $(g',h'), (g'',h'') \in V(G\circ H)$ are adjacent
 if either $g'g''\in E(G)$ or $g'=g''$ and
 $h'h''\in E(H)$.

\ 

Two less commonly known, yet still significant graph parameters, are defined as follows:
\begin{itemize}
\item $a(G)$ represents the \textit{vertex arboricity} of a graph $G$, which is defined as the minimum number $a$ of vertex classes in a partition $(V_1,\dots,V_a)$ of $V(G)$ such that each $V_i$ induces a forest in $G$.

\item $d(G)$ denotes the \textit{degeneracy} of a graph $G$. It is the smallest non-negative integer $d$ such that every induced subgraph $G'$ of $G$ satisfies $\delta(G')\leq d$.
\end{itemize}

\subsection{Our results}

In Section 2, we study how basic graph operations---edge or vertex deletion, union and vertex-disjoint union, lexicographic product, or taking the line graph---act on the robust chromatic number.

Then in Section 3, we compare $\chi_1(G)$ to some of the most commonly considered graph parameters. We summarize our findings in the theorem below.

\bthm\label{param}
\

\begin{enumerate}
\item 
For every graph $G$ we have
$$\left\lceil \frac{\chi(G)}{3} \right\rceil \le \cir(G) \le \chi(G)$$
and
$$\left\lceil \frac{\omega(G)}{3} \right\rceil \le \omi(G) \le \omega(G).$$
All these bounds are tight, for all possible values of $\chi$ and $\omega$.
\item
For every isolate-free graph $G$,
$$\theta(G)\leq\theta _1(G)\leq 3\theta(G)$$
and the upper bound is tight.
\item
For every graph $G$ the bounds $$a(G)/2 \leq \cir(G) \leq a(G)$$
 are valid and tight. 
\item
Let $k$ be any positive integer.
If $\Delta(G) < 3k$, then
$$\cir(G)\leq k,$$
that is,
$$\cir(G) \leq \left\lceil \frac{\Delta(G)+1}{3} \right\rceil.$$
Moreover, the bounds are tight for both $\Delta$ and $\cir$
as there exist graphs $G_k$ with $\Delta(G_k)=3k$ and $\cir(G_k)=k+1$.
\item
Every $d$-degenerate graph $G$ has $$\cir(G) \leq d/2 + 1.$$
Moreover, this upper bound is tight as for every $k\geq 1$ there exists a graph $H_k$
such that $H_k$ is $2k$-degenerate and $\cir(H_k)=k+1$.
\item
If $\Delta(G)>1$, then
$$\cir'(G)\leq\chi'(G)-2.$$
Moreover,
$$\delta(G) - 2 \leq \cir'(G) \leq \Delta(G) - 1.$$
All these bounds are tight.
\end{enumerate}
\ethm
Separate points of Theorem \ref{param} will be proved in different subsections of Section~3.

\medskip 

In the present context it is natural to introduce the following two algorithmic problems:

\begin{quotation}

{\sc Robust $k$-colorability}

\textbf{Input:} Graph $G=(V,E)$, natural number $k$.

\textbf{Question:} Is $\cir(G) \leq k$\,?

\end{quotation}

\begin{quotation}
{\sc Robust coloring}

\textbf{Input:} Graph $G=(V,E)$.

\textbf{Solution:} The value of $\cir(G)$.
\end{quotation}

Due to the next result, which we prove in Section 4, it would be of great interest to identify
 graph classes in which $\cir$ can be determined efficiently.

\bthm\label{complexity}
For every natural number\/ $k\geq 3$, the\/
 {\sc Robust $k$-colorability} problem is\/ \npc.
Moreover, {\sc Robust coloring} is not approximable within\/
 $O(|V|^{1/2-\ve})$ for any real\/ $\ve > 0$, unless\/ {\sf P $=$ NP}.
\ethm

%

On the positive side, we prove that $\ali,\omi,\cir,\thi$ are computable on graphs of bounded treewidth in linear time.

\msk

Many further questions are raised in the concluding section.

\paragraph{Vertex partition vs.\ edge decomposition.}

We close the introduction with observations comparing the
 quasi-unicyclic partitions of vertex sets and edge sets.
The parameter $\cir$ asks about the minimum number of classes
 in a vertex partition of a graph into sets that \emph{induce}
  quasi-unicyclic subgraphs.
It turns out that \emph{edge decompositions into subgraphs} of this kind have a completely different nature and can be handled in a more efficient way.

\bthm
For any graph\/ $G=(V,E)$ the minimum number of quasi-unicyclic
 subgraphs\/ $H_1,\dots,H_k$ with\/
  $E(H_1)\cup\cdots\cup E(H_k) = E$
 is equal to
  $$
    \max_{U\subseteq V} \left\lceil \frac{e(G[U])}{|U|} \right\rceil,
  $$
 where\/ $e(G[U])$ denotes the number of edges induced by\/ $U$ in\/ $G$.
Moreover, the minimum\/ $k$ and a corresponding edge decomposition
 can be determined in polynomial time.
\ethm

\pf
It follows from a result of Hakimi \cite[Theorem 4]{H} that
 if $G=(V,E)$ is an undirected graph and $t$ is a positive integer
 such that each $Y\subseteq V$ induces at most $t\cdot|Y|$ edges
 in $G$, then $G$ admits an orientation with maximum out-degree
 at most $t$.
For later history, references and short proofs of the general
 form of Hakimi's theorem we refer to Section 2.2 of \cite{STV}.
Once an orientation of this kind with
 $t=\max_{U\subseteq V} \left\lceil \frac{e(G[U])}{|U|} \right\rceil$
 is at hand for $G$, we can partition the edge set into
 $t$ classes so that the out-going edges at each vertex belong
 to mutually distinct classes.
Then each class forms a graph with all vertices having out-degree
 0 or~1, thus each edge class is a quasi-unicyclic graph.

Concerning the cited known results it is also known that
 an orientation minimizing the maximum out-degree
 can be obtained by using maximum matching algorithms in bipartite graphs.
This yields a solution in polynomial time.
\qed

\section{Elementary graph operations}

The following result collects the effect of some frequently
 studied operations on graphs.

\bthm   \label{t:oper}
\begin{enumerate}[label=(\roman*)]
\item 
\textrm{[Edge deletion, vertex deletion.]} 

The invariants $\ali, \omi, \cir, \thi$ are monotone with respect to graph inclusion, in the following way.

If $G$ is any graph and $H$ is a \underline{spanning} subgraph of $G$, then
$$
   \ali(G) \leq \ali(H),\qquad
    \thi(G) \leq \thi(H),
 $$

 and if $H$ is any subgraph of $G$, then
 $$
   \omi(G) \geq \omi(H),\qquad
    \cir(G) \geq \cir(H).
 $$

\smallskip

\item 
\textrm{[Vertex-disjoint union]}

If $G$ is disconnected and $G_1, \dots, G_k$ are its connected components, then
 $$
   \cir(G) = \max_{1\leq i\leq k} \cir(G_i) ,\qquad
   \omi(G) = \max_{1\leq i\leq k} \omi(G_i) ,
 $$
 $$
   \ali(G) = \sum_{i=1}^k \ali(G_i) ,\qquad
   \thi(G) = \sum_{i=1}^k \thi(G_i) .
 $$

\smallskip

\item 
\textrm{[Union of two graphs]}

If $V(G_1)=V(G_2)$, then
 $$
   \cir(G_1\cup G_2) \leq
    \min \left\{ \chi(G_1)  \cir(G_2) , \,
     \cir(G_1)  \chi(G_2) \right\} ,
 $$
 
 and the bound is tight.

\smallskip

\item 
\textrm{[General graph union]}

If $G_1, \dots, G_k$ are graphs on the same vertex set, then
 $$
   \cir(G_1\cup\cdots\cup G_k) \leq
    (2k+1) \prod_{i=1}^{k} \cir(G_i) .
 $$
Moreover, there exist graphs $G_1, \dots, G_k$ such that
 $$
   \cir(G_1\cup\cdots\cup G_k) \geq
    \frac{2k+1}{3} \prod_{i=1}^{k} \cir(G_i) .
 $$

\smallskip

\item 
\textrm{[Lexicographic product]}

For any two graphs $G$ and $H$ we have
 $\cir(G\circ H) \leq \chi(G) \cdot \cir(H)$.

\smallskip

\item 
\textrm{[Line graph]}

Contrary to the notion of proper coloring, the unary operation $G\mapsto L(G)$ does not admit a direct correspondence between robust edge colorings of $G$ and robust vertex colorings of $L(G)$;
 in particular, $\cir'(G)=\cir(L(G))$ does not hold in general.
\end{enumerate}
\ethm

\pf
Proof of $(i)$: it suffices to note that the inequalities
 $\alpha(G) \leq \alpha(G-e)$,
 $\chi(G) \geq \chi(G-e)$,
 $\omega(G) \geq \omega(G-e)$,
 $\theta(G) \leq \theta(G-e)$
  hold for every graph $G=(V,E)$ and every edge $e\in E$.

\msk

Proof of $(ii)$: the equalities follow by choosing an optimal 1-selection in each component of $G$.

\msk

Proof of $(iii)$: it suffices to take a 1-selection in one of the two graphs,
and apply the fact that $\chi$ is a submultiplicative function. 
 Tightness is shown by the following example.
 For $p,q\ge 3$ with $3p\geq q$, we partition $G=K_{3pq}$ into the complete $q$-partite graph $G_1=K_{3p,3p,\dots,3p}$ and the disjoint union of $q$ cliques, $G_2=qK_{3p}$. Then $\chi_1(G)=pq$, $\chi(G_1)=\chi_1(G_1)=q$, $3\chi_1(G_2)=\chi(G_2)=3p$, and so $pq=\min\{\chi(G_1)\chi_1(G_2),\chi_1(G_1)\chi(G_2)\}=\min\{qp,q\cdot 3p\}$.

\msk

Proof of $(iv)$:
For $i=1,\dots,k$ let $F_i$ be the edge set in an optimal
 1-selection $f_i$ on $G_i$
 (optimal in the sense that $\chi(G_i-F_i)=\cir(G_i)$).
Then any subset $X$ of the vertex set induces at most $|X|$ edges
 of $F_i$, hence the average degree in every induced subgraph of
  $$
    F := F_1\cup\cdots\cup F_k
  $$
 is at most $2k$.
As a consequence, $\chi(F) \leq 2k+1$.
Thus, writing the union in the form
 $$
   G_1\cup\cdots\cup G_k = F \cup (G_1-F_1) \cup \cdots \cup (G_k-F_k)
 $$
  and applying the multiplicative property of $\chi$,
  the upper bound follows.

A simple example showing the claimed lower bound is
 obtained from a Hamiltonian decomposition of the complete graph $K_{2k+1}$.
This yields $k$ Hamiltonian cycles, which
 we can take as $G_1,\dots, G_k$.
Then of course $\cir(G_i)=1$ holds for every $i$,
 while, as Example \ref{ex:kn} shows, we have $\chi_1(K_{2k+1})=\lceil \frac{2k+1}{3}\rceil$.

\msk

Proof of $(v)$:
For every $g\in V(G)$ the vertex subset
 $V_g := \{ (g,h) \mid h\in V(H) \}$
 induces a subgraph isomorphic to $H$ in $G\circ H$.
Choose an optimal 1-selection $f:V(H)\to E(H)$ of $H$, and
 copy it into each $V_g$, hence creating a 1-selection
 $f_{G\circ H}$.
Removing the edge set $f_{G\circ H} V(G\circ H)$
 from $G\circ H$, each $V_g$ admits a partition into
 independent sets $V_{g,i}$ for $i=1,\dots,\cir(H)$.
Hence we can consider a proper vertex coloring
 $\vp:V(G)\to\{1,\dots,\chi(G)\}$ of $G$ with the
  minimum number of colors, and decompose
 $V(G\circ H)$ into the sets $\bigcup_{g\in S} V_{g,i}$
 where $S$ runs over the color classes of $\vp$.
Clearly, each of those $\chi(G) \cdot \cir(H)$ sets is
 independent in $(G\circ H) - f_{G\circ H} V(G\circ H)$,
  verifying the validity of the assertion.

\msk

Proof of $(vi)$:
As a small example, $L(K_4)=K_6-3K_2$, hence $\cir(L(K_4))=2$,
 while removing the 1-selection $C_4\subset K_4$ we obtain
  $2K_2$, therefore $\cir'(K_4)=1$.
As an infinite class, the stars $K_{1,m}$ have
 $\cir'(K_{1,m})=1$ and
  $\cir(L(K_{1,m}))=\cir(K_m)=\lceil m/3 \rceil$.
\qed

\bsk

It should be noted that the two invariants in $(v)$ are
 not interchangeable; that is, $\cir(G) \cdot \chi(H)$ is
 not an upper bound on $\cir(G\circ H)$.
Simple counterexamples are obtained by taking
 sufficiently large edgeless graphs for $H$.

\section{Comparison with other graph invariants}

In this section, we compare the robust parameters $\omi$, $\thi$,
 and $\cir$ with several important graph invariants.

\subsection{Clique number}

\begin{prop} \label{f:chi}
For every graph\/ $G$ we have
\begin{equation} \label{eq:chi}
 \left\lceil \frac{\chi(G)}{3} \right\rceil \le \cir(G) \le \chi(G)
\end{equation}
 and
\begin{equation} \label{eq:omega}
 \left\lceil \frac{\omega(G)}{3} \right\rceil \le \omi(G) \le \omega(G).
\end{equation}
All these bounds are tight, for all possible values of\/
 $\chi$ and\/ $\omega$.
\end{prop}
\pf 
The upper bounds follow directly from the definitions.
Also the lower bound in (\ref{eq:omega}) is implied by the
 property of complete graphs shown in Example~\ref{ex:kn}.
For the lower bound in (\ref{eq:chi}), let $k=\cir(G)$ and
 consider the color classes $V_1,\dots,V_k$ of a subgraph $G_f$
  of $G$, where $k=\chi(G_f)=\cir(G)$. 
By Proposition~\ref{p:basic}$(i)$, the induced subgraph $G[V_i]$
 is unicyclic and therefore $3$-colorable for all $i \in[k]$. 
Thus, each $V_i$ can be partitioned into at most three sets that are independent in $G$, and this way we obtain a color partition of $G$ with at most $3k=3\cir(G)$ color classes. 
This shows $\chi(G)\le 3\cir(G)$ and the lower bound follows.

Tightness is shown for any value of $\chi$ by
 Examples \ref{ex:kn} and \ref{ex:kttt}.
\qed
\bsk

In fact, there are much wider classes of graphs establishing
 equality in either side of (\ref{eq:chi}), as it will be proved
  in the sequel.

\subsection{Clique covering number}

Here we deal with the robustness parameter $\thi(G)$
 corresponding to $\theta (G)$.

\begin{prop}   \label{p:theta}
For every isolate-free graph\/ $G$,
$$\theta(G)\leq\theta _1(G)\leq 3\theta(G)$$
and the upper bound is tight.
\end{prop}
\pf
It is evident to note that for every graph $G$, the inequality $\thi(G)\geq \theta(G)$ holds.
For the proof of the upper bound $3\theta(G)$ let
 us take a minimum cover (co-coloring) of $G$ and let $Q$ be an arbitrary class in it.
 Consider any $1$-selection $S$ of $G$ as a graph. It is more convenient to work with the complement of the graph. The subgraph of $S$ induced by $Q$ has chromatic number at most $3$ as all the cycles and forests are $3$-colorable. Thus, if we return to $G$ and we delete $S$ from $Q$, the subgraph obtained has co-chromatic number at most $3$. 
 Applying this for any class, we obtain the claimed bound.
 
Tightness is shown by the vertex-disjoint union of any cliques of sizes at least 3.
Indeed, a triangle can be removed from each of those cliques.
Then every original clique will need at least three cliques in a clique cover after the removal of a suitable 1-selection.
\qed

\subsubsection{An open problem}

We formalize here the following conjecture.
\begin{conj}   \label{cj:theta}
The lower bound\/ $\thi \geq \theta$ is not tight, except for edgeless graphs.
\end{conj}

In order to provide partial results in this direction,
in the rest of this subsection we analyze the properties of a hypothetic counterexample.
We shall use the following term for it.

\begin{defi}
We call a graph\/ $G$ \emph{exact} if\/
$\theta _1(G) = \theta (G)$.
\end{defi}

\subsubsection{Properties of exact graphs}

Except where explicitly stated otherwise, $G$ is assumed to be an exact graph throughout this sub-subsection.


%
\begin{defi} \label{crdef}
A graph\/ $G$ is\/ \emph{critically\/ $k$-co-chromatic} if\/
 $\theta (G)=k$ but\/ $\theta (G-x)=k-1$ holds for every\/ $x\in V(G)$.
Further,\/ $G$ is\/ \emph{critical} if it is critically\/
 $k$-co-chromatic for some\/ $k$.
\end{defi}

\begin{prop}\label{notcr}
$G$ is critical.
\end{prop}

\pf
If $G$ is not critical, there exists  a vertex $x$ such that $\theta (G-x)=\theta (G)$. 
Let us denote by $X$ the set of all edges incident to $x$. In the graph $H:=G-X$ any clique cover contains $\{x\}$ and so $\theta (H)=\theta (G-x)+1$. The edge set $X$ can be extended to a spanning forest $F$ and 
$$\theta (G-F)\geq \theta (G-X)=\theta (H)=\theta (G-x)+1=\theta (G)+1.$$
Here $G-F$ is a removed graph, thus $\theta _1(G)\geq \theta (G-F)>\theta (G)$.
\qed

\begin{prop} \label{perfnotcr}
 $\theta (G) >\alpha (G)$.
\end{prop}

\pf
Suppose by way of contradiction that $\theta (G) = \alpha (G)$ holds.  Taking $x\in V(G)$ arbitrarily, and using both equality and criticality,
  $$\alpha (G-x) \leq \theta (G-x)\leq \theta (G)-1=\alpha (G)-1$$
Consequently, $x$ is contained in every maximum independent set of the graph. As $x$ is arbitrary, this implies that $G$ is edgeless, a contradiction.
\qed \bsk

Next, we present a surprising fact.

\begin{cor}  \label{perfdone}
$G$ is imperfect.
\end{cor}
\pf
By way of contradiction, assume $G$ is perfect. Then, from the Weak Perfect Graph Theorem    \cite{L1}, we have $\theta(G)=\alpha(G)$, contradicting Proposition \ref{perfnotcr}.
\qed
\begin{cor}  \label{P4perf}
There exists an induced\/ $P_4$ in\/ $G$.
\end{cor}
\pf
A very old result of Seinsche \cite{Sei} states that if a graph has no induced $P_4$ then it is perfect. So, we obtain the result from Corollary \ref{perfdone}.
\qed

\begin{prop} \label{abc}
The following is impossible for\/ $G$:  For some vertex triple\/ $a,b,c$ of\/ $G$, in any
 co-coloring of\/ $G-a$, the vertices\/ $b$ and\/ $c$ have the same co-color.
\end{prop}
\pf
Suppose $G$ is exact and consider  the edge set $X$ consisting of the edge $bc$ together with all edges incident to $a$.  We extend $X$ to  a $1$-selection $S$. Thus, in an arbitrary minimum co-coloring of $G-S$, $\{a\}$ will yield a singleton class, moreover, $b$ and $c$ cannot have the same co-color, a contradiction.
\qed

\bsk

From Proposition \ref{abc} we obtain:

\begin{cor} \label{nique}
For any vertex\/ $a$ of\/ $G$, the  graph\/ $G-a$ is not uniquely co-colorable.
\end{cor}
 
\begin{defi}
A graph is \emph{partitionable} if for every vertex\/ $a$, the graph\/ $G-a$ can be represented as a rectangle where the rows are stable sets and the columns are cliques. 
\end{defi}

It can be proved (see \cite{L2}) that in a partitionable graph $G$, the size of the rows is $\alpha (G)$, and the size of the columns is $\omega (G)$.

\begin{cor} \label{partcol}
$G$ is not partitionable.
\end{cor}
\pf
It is a well-known (but nontrivial) fact proved in \cite{Padb} that in a partitionable graph, $G-a$ is uniquely co-colorable for any $a$.
So, we can apply Corollary \ref{nique}.
\qed

\bigskip

From the results above it follows that chordless cycles and their complements
 are non-exact.
(For even order, they are perfect, for odd order they are partitionable.
Certainly, non-exactness can also be verified by
 elementary direct proofs.)



Since $G$ is critical due to Proposition \ref{notcr}, its complementary graph is chromatic critical, and the following statement is true.
    

\begin{claim}  \label{coconn}
Let\/ $G$ be a counterexample for Conjecture \ref{cj:theta}, having minimum number of vertices. Then both\/ $G$ and its complement are connected.
\end{claim}

%
\bthm
Every edge of\/ $G$ is contained in at least two triangles.
\ethm
\pf
Let $\theta (G)=k$.
    By way of contradiction, suppose that an edge $x_1x_2=e\in E(G)$ occurs in at most one triangle. Then the subgraph $S$ with the edges containing $x_1$ or $x_2$ (or both) is quasi-unicyclic. Since $G$ is exact, $\theta(G-S)=\theta (G)$ holds; that is, $G-S$  has some co-coloring $\cal C$ with $k$ colors. According to the construction of $S$, $x_1$ and $x_2$ are isolated in $G-S$, consequently they necessarily yield one-element co-color classes in $\cal C$. Moreover, within $G[V-\{x_1, x_2\}]$,  $\cal C$ has $k-2$ classes that form cliques in $G$, too. But we may attach the clique consisting of $x_1$ and $x_2$, thus  yielding a ($k-1$)-co-coloring in $G$, a contradiction. 
    
    \qed
    
\begin{prop}\label{theta4}
$\theta (G)\geq 4$.
\end{prop}

\pf
Otherwise, $\theta (G)\leq 3$.
If $\theta (G)$ is at most $2$, then $\overline{G}$ is bipartite, and consequently $G$ is perfect, which is impossible. This implies that $\theta (G)=3$ and $\overline{G}$ is critically $3$-chromatic. Then, clearly,  $\overline{G}$  is a chordless odd cycle, which we already excluded by Corollary \ref{partcol}.
\qed

\begin{defi}
A vertex set\/ $W$ is an\/ \emph{inducing set} if it induces a quasi-unicyclic subgraph in\/ $G$. We denote
$$\iota (G):= \max \{|W| : W \ \textrm{is an inducing set in} \ G\} .$$
\end{defi}
\begin{prop}    \label{jota}
$\theta (G) \geq \iota (G)$.
\end{prop}
\pf
Let us take an inducing vertex set $W$ in $G$ 
with $|W|=\iota (G)$.
The set of edges induced by $W$ forms a quasi-unicyclic subgraph, and it can be extended to a $1$-selection $f$ of $G$. If we delete all the edges of $f$ from $G$, we obtain a $1$-removed graph $G_f$. Hence $W$ will be independent in $G_f$. Consequently,
$$\theta (G) = \theta _1 (G) \geq \theta (G_f) \geq \alpha (G_f)\geq |W|,$$
proving the assertion.
\qed

\begin{rmk}
As an illustration, let\/ $\theta (G)=4$. In this case, no inducing set of size\/ $5$ can exist in\/ $G$.
For the collection\/ $\cal C$ of graphs to be thus forbidden as an induced subgraph, we give a list of  graphs of order\/ $5$, maximal in\/ $\cal C$  with respect to their edge sets: 
\msk

$C_5$,  $C_4$ plus one leaf, and three graphs constructed from a\/ $C_3$: 

\msk

The bull (two leaves attached on different vertices), two leaves attached

on the same vertex, and finally, an attached path of length\/~$2$.
\msk

\noindent
It follows that any spanning subgraph of any member of this list is forbidden in\/~$G$.
\end{rmk}

\begin{defi}
A vertex set\/ $D$ is\/ \emph{dominating} in\/ $G$ if for every\/ $x\in V-D$,\/ $x$ has at least one neighbor in\/ $D$.
\end{defi}

\begin{defi}
An edge\/ $ab$ is\/ \emph{dominating} if the set\/ $\{a,b\}$ is dominating.
\end{defi}

\begin{prop} \label{indd}
If\/ $\theta (G)\leq 4$ then\/ $G$ contains a dominating edge.
\end{prop}

We prove a stronger statement.
Recall from Corollary \ref{P4perf} that $G$ contains at least one induced $P_4$.
The next assertion implies that the middle edge of any induced $P_4$ in $G$ is a dominating edge.
%
%

\begin{prop} 
If\/ $\theta (G)\leq 4$ then for each vertex set\/ $P$ which induces a\/ $P_4$ and for each\/ $x\notin P$,\/ $x$ has at most one non-neighbor in\/ $P$.
\end{prop}
\pf
If a vertex $x$ has at least two non-neighbors in $P$, then $P\cup \{x\}$ induces a quasi-unicyclic subgraph of order $5$.
By Proposition \ref{jota} this would imply $\theta(G)\geq \iota(G) \geq 5$, a contradiction.
\qed

\subsection{Vertex arboricity}

\begin{prop}   \label{p:arb}
 For every graph\/ $G$ the bounds\/ $a(G)/2 \leq \cir(G) \leq a(G)$
 are valid and tight.
\end{prop}

\pf
 The upper bound follows directly from the definitions.
Moreover it is tight because $\cir(G) = 1 = a(G)$ holds whenever
 $G$ is a tree.
For arbitrarily large values of $a(G)$, we refer to
 Example \ref{ex:kttt}: a complete multipartite graph with any
 number $k$ of vertex classes and  more than $k$ vertices in
 each class satisfies the equality $\cir=\chi=k$, hence its
 vertex arboricity is also the same.

For the lower bound we use Proposition \ref{p:basic}$(i)$ and
 observe that the vertex set of each omitted cycle under an
 optimal edge-selecting function $f$ can be partitioned into
 two paths, hence obtaining a coloring of $G$ such that each
 color class induces a tree.
Tightness is shown e.g.\ by any graph $G$ in which each connected
 component is a cycle.
Then we have $\cir(G)=1$ and $a(G)=2$.

We can give constructions $G_k$ for general $\cir=k$, too.
Let $V_1,\dots, V_k$ be mutually disjoint sets of size $3k$ each.
Put a complete bipartite graph $K_{3k,3k}$ between any two
 $V_i,V_j$ and put $k$ disjoint triangles inside each $V_i$.
We prove that this graph has $a(G_k)=2k=2\cir(G_k)$.

Suppose that $X_1\cup\cdots\cup X_a = V(G_k)$ is a vertex
 partition into $a=a(G_k)$ classes, such that each $X_\ell$
 induces an acyclic subgraph.
For a distinction, we call a $V_i$ a part of $G$, and an
 $X_\ell$ a class of $G$.
No class can meet more than two parts, otherwise, it would
 induce at least one triangle.
Hence we may have ``single classes'' entirely contained in a
 part, and ``double classes'' that meet two parts.
We first consider the single classes.

If a part contains more than one single class, then we may
 assume without loss of generality that it is the union of
 exactly two single classes.
We remove all those parts and classes, say $k''$ parts and
 $2k''$ classes.
The remaining graph, say $G'$, has $k':=k-k''$ parts and
 vertex arboricity $a':=a-2k''$.
We need to prove that $a'\geq 2k'$ holds.
Let $a'=s+d$, where $s$ and $d$ denote the number of single
 classes and the number of double classes, respectively.
The size of a single class is at most $2k$, and there can be
 at most $k'$ of them; while the size of a double class is
 at most $k+1$, because it can meet one of the two parts in
 just one vertex (in order to avoid a $C_4$) and can contain
 at most one vertex from each triangle of the other part.
Since all of the $3kk'$ vertices must be covered, we obtain
 $$
   a' = s+d \geq s+ \frac{3kk'-2ks}{k+1} =
    \frac{3kk'-(k-1)s}{k+1} \geq k' \cdot \frac{2k+1}{k+1} =
     2k' - \frac{k'}{k+1},
 $$
that means $a'\geq 2k'$ as $a'$ is an integer.
\qed

\subsection{Vertex degrees}

\bthm [Maximum degree]
Let\/ $k$ be any positive integer.
If\/ $\Delta(G) < 3k$, then\/ $\cir(G)\leq k$; that is,\/
 $\cir(G) \leq \left\lceil \frac{\Delta(G)+1}{3} \right\rceil$.
Moreover, the bounds are tight for both\/ $\Delta$ and\/ $\cir$
 as there exist graphs\/ $G_k$ with\/
 $\Delta(G_k)=3k$ and\/ $\cir(G_k)=k+1$.
\ethm

\pf
Beginning with the assertion on tightness, the
 complete graphs $G_k=K_{3k+1}$ are suitable examples.

For the assertion on $\cir$, let $G=(V,E)$ be a graph
 with maximum degree at most $3k-1$.
We take a vertex partition $(V_1,\dots,V_k)$ of $G$
 such that the total number of edges joining distinct
 classes $V_i,V_j$ ($1\leq i<j\leq k$) is as large
 as possible.
Then each class induces a subgraph of maximum degree
 at most 2.
Indeed, if $v\in V_i$ has at least three neighbors
 inside $V_i$, then at most $3k-4$ edges join $v$
 to the other $k-1$ classes, hence there is a class
 $V_j$ in which $v$ has at most two neighbors.
Re-defining then $V_i:=V_i\setminus\{v\}$ and
 $V_j:=V_j\cup\{v\}$ we obtain a partition with
 more crossing edges, a contradiction.
It follows that each class induces a union of paths
 and cycles, therefore a 1-selection can contain all
 edges inside the $k$ classes, thus $\cir(G)\leq k$.
\qed

\bthm [Degeneracy]   \label{t:dege}
Every\/ $d$-degenerate graph\/ $G$ has\/ $\cir(G) \leq d/2 + 1$.
Moreover, this upper bound is tight as
 for every\/ $k\geq 1$ there exists a graph\/ $H_k$
 such that\/ $H_k$ is\/ $2k$-degenerate and\/ $\cir(H_k)=k+1$.
\ethm

\pf
Consider a graph $G$ with degeneracy number $d$. Let $v_1,v_2,\dots,v_n$ be an enumeration of the vertices of $G$ such that every $v_i$ has at most $d$ neighbors in $\{v_j:j<i\}$. We define a 1-selection $f:V(G)\rightarrow E(G)$ and a coloring $c:V(G)\rightarrow \{1,2,\dots, \lfloor d/2\rfloor+1\}$ simultaneously. We let $c(v_1)=1$ and $f(v_1)$ an arbitrary edge of $G$ incident to $v_1$. Suppose we have defined $c$ and $f$ for all vertices $v_1,\dots,v_i$ such that $c$ is a proper coloring of $G[v_1,\dots,v_i]_f$. Then as $v_{i+1}$ has at most $d$ neighbors in $v_1,\dots,v_i$, there must exist a color class $c^{-1}(j)$ for some $j\le \lfloor d/2\rfloor +1$ such that $v_{i+1}$ has at most one neighbor in $c^{-1}(j)$. We then let $c(v_{i+1})=j$ and define $f(v_{i+1})$ to be the edge joining $v_{i+1}$ to its only neighbor in $c^{-1}(j)$ (if it exists, otherwise $f(v_{i+1})$ can be an arbitrary edge incident to $v_{i+1}$). Clearly, once $c$ and $f$ are defined on the entire graph, $c$ is a proper coloring of $G_f$. This finishes the proof of the upper bound $d/2+1$.

Tightness for $d=2$, that is $k=1$, is clear by $H_1:=K_4-e$.
A general construction will have vertex set
 $V=V(H_k)=V_0\cup V_1\cup\cdots\cup V_{k'}$
  where $k'=\lceil (2k+1)/3 \rceil$ will suffice.
The subgraph induced by $V_0$ is $K_{2k}$.
For each $1\leq i\leq k'$ the set $V_i$ is independent and has size
 $(k+1)\cdot\binom{|V_0|+\ldots+|V_{i-1}|}{2k}$.
Each $v\in V_i$ has $2k$ neighbors in
 $\bigcup_{j=0}^{i-1} V_j$, and any $2k$ vertices of
 $\bigcup_{j=0}^{i-1} V_j$ have $k+1$ common neighbors in $V_i$.
This graph $H_k$ clearly has degeneracy number $2k$, and so the
 first part of the theorem guarantees $\cir(H_k)\leq k+1$.

Suppose for a contradiction that $\cir(H_k)\leq k$, and let
 $V=X_1\cup\cdots\cup X_k$ be a vertex partition where
 each $X_j$ induces a quasi-unicyclic graph. Observe that $K_4^-:=K_4-e$ is not quasi-unicyclic. For an $i\ge 0$ let us write $c_i=|\{j:X_j\cap (\cup_{h=0}^iV_h)\neq \emptyset\}|$ and $p_i=|\{j:H_k[X_j\cap (\cup_{h=0}^iV_h)]$ contains an edge$\}|$. As $K_k[V_0]$ is complete, we have $c_0+p_0\ge \lceil \frac{4k}{3}\rfloor$ and $c_0\ge \lceil \frac{2k}{3}\rceil$. We claim that as long as $p_i<k$, we have $c_i+p_i<c_{i+1}+p_{i+1}$. Indeed, consider a set $D\subseteq \cup_{h=0}^iV_h$ that contains an edge in all $p_i$ possible colors and a vertex from all possible $c_i$ color classes. There exist a set $N$ of $k+1$ vertices in $V_{i+1}$ that are joined to all vertices of $D$. As $K_4^-$ is not quasi-unicyclic, $N$ can contain at most one vertex from each color class with an edge in $D$. As $p_i<k$, there exists a vertex $x\in N$ that is not of these colors. If its color class is completely new, then $c_{i+1}>c_i$ increases; and if it appears before, so in $D$, then $p_{i+1}>p_i$. As $c_i\le k$ for all $i$, $c_0+p_0\ge \lceil \frac{4k}{3}\rfloor$ and $c_0\ge \lceil \frac{2k}{3}\rceil$ imply $p_i=k$ for some $i\le \frac{2k}{3}$. 
 
 Finally, we claim that if $p_i=k$, then the color classes $X_1,X_2,\dots,X_k$ cannot be extended to $V_{i+1}$. To see this, consider again a set $D'$ of $2k$ vertices that contains an edge from each color in $\cup_{h=0}^iV_i$, and let $N'$ be its joined neighborhood of $k+1$ vertices in $V_{i+1}$. By the pigeon-hole principle there exist two vertices $x,y$ in the same color class, say in $X_1$. Then together with the edge $e$ in $D\cap X_1$, they form a $K_4^-$ in $X_1$, contradicting the fact that $H_k[X_1]$ is quasi-unicyclic.
\qed

\subsubsection{Consequences for planar and outerplanar graphs}

In this extremely short subsection, we derive two consequences on planar graphs, whose coloring properties are among the most classical issues in graph theory.

\bthm   \label{t:pla}
$(i)$~~If\/ $G$ is an outerplanar graph, then\/
 $\cir(G)\leq 2$.

$(ii)$~~If\/ $G$ is a planar graph, then\/
 $\cir(G)\leq 3$.
\ethm

\pf
Both parts are consequences of Theorem \ref{t:dege}.
Every outerplanar graph is 2-degenerate, hence $(i)$
 follows by taking $d(G)=2$.
Moreover, every planar graph is 5-degenerate, hence $(ii)$
 follows by taking $d(G)=5$.
\qed

\subsection{Chromatic index}

\bthm
If\/ $\Delta(G)>1$, then $\cir'(G)\leq\chi'(G)-2$. Moreover,
 $$
   \delta(G) - 2 \leq \cir'(G) \leq \Delta(G) - 1 .
 $$
All these bounds are tight.
\ethm

\pf
To prove the upper bound $\chi'(G)-2$ we consider an edge coloring
 $\psi$ with $\chi'(G)$ colors.
Choose two color classes, say $E_1$ and $E_2$.
Then in $E_1\cup E_2$, each connected component is a path or a
 cycle, hence $E_1\cup E_2$ can be made a 1-selection $f$.
The restriction of $\psi$ to $E\smin (E_1\cup E_2)$ properly
 edge-colors $G_f$ with $\chi'(G)-2$ colors.
This also implies $\cir'(G) \leq \Delta(G)-1$ by Vizing's theorem.

For the lower bound $\delta(G) - 2$ we observe that removing
 at most $|V(G)|$ edges makes the vertex degrees decrease by
 at most 2 on average.
Thus, there remains a vertex with degree of at least $\delta(G) - 2$,
 implying $\chi'(G_f) \geq \delta(G)-2$ for every 1-selection $f$.

Regular graphs of type 1 have $\delta(G)=\Delta(G)=\chi'(G)$,
 and every color class in an optimal edge coloring is a
 perfect matching.
Hence the removal of two color classes decreases all of these
 parameters with exactly 2.
Tightness of $\cir'(G) \leq \Delta(G) - 1$ is shown e.g.\ by
 complete graphs of odd order.
\qed

\bsk

On the other hand, it has to be noted that there is no lower bound
 on $\cir'(G)$ in terms of $\Delta(G)$.
This fact is shown by trees of any large maximum degree, which have
 $\cir'=0$.

\section{Algorithmic complexity}   \label{s:algo}

In the first part of this section we prove that it is hard to compute, and even to
 approximate, the robust chromatic number of a generic input graph.
After that, we show how all the four parameters $\ali,\omi,\cir,\thi$ are 
computable in linear time on graphs of bounded treewidth.

For the \nph ness result, we restate Theorem \ref{complexity}:

\begin{itemize}
    \item[ ] \emph{For every natural number\/ $k\geq 3$, the\/
 {\sc Robust $k$-colorability} problem is\/ \npc.
Moreover, {\sc Robust coloring} is not approximable within\/
 $O(|V|^{1/2-\ve})$ for any real\/ $\ve > 0$, unless\/ {\sf P $=$ NP}.}
\end{itemize}

\pf
We begin with the observation that {\sc Robust $k$-colorability}
 is in the class \np.
A certificate, that can be verified in polynomial time, is a
 vertex $k$-partition such that each class induces a
 quasi-unicyclic graph.
Here $k$ is not required to be fixed,
 it may also depend on the order of the input graph.

To prove the hardness results, we apply reduction from the
 corresponding problems on proper vertex colorings of graphs.
As it is well known, for every $k \geq 3$ it is \npc\ to decide
 whether a generic input graph is $k$-colorable.
Now, for any $G=(V,E)$ of order $n$, we substitute each vertex $v$
 of $G$ with an independent set $S_v$ of size $n+1$;
  if two vertices $v,w\in V$ are adjacent, the edge $vw$ is enlarged
  to $K_{n+1,n+1}$, otherwise no edges are drawn between the
  corresponding two $(n+1)$-sets.
In this way a graph $G^+$ of order $n(n+1)$ is obtained, and the
 transformation takes polynomial time.
We claim:
 $$
   \cir(G^+) = \chi(G^+) = \chi(G) .
 $$

The second equality is straightforward since $G$ is a subgraph of $G^+$, and
 on the other hand, every proper coloring of $G$ can be enlarged in
 a natural way to a proper vertex coloring of $G^+$ with the same
 number of colors.

To verify the first equality we observe that picking one vertex from
 each set $S_v$ in all possible ways, we obtain $(n+1)^n$ distinct
 subgraphs isomorphic to $G$.
Each edge is contained in $(n+1)^{n-2}$ of those subgraphs.
Hence removing a 1-selection $f(V(G^+))$ from $G^+$ we can destroy
 no more than $n(n+1)^{n-2}$ copies of $G$, consequently we still
 have $G \subset (G^+ - f(V(G^+))$.
Thus, $\cir(G^+) \geq \chi(G) = \chi(G^+) \geq \cir(G^+)$.
This implies equality and finishes the proof of \npc ness.

To prove inapproximability, we cite Zuckerman's important result
 \cite{Z} stating that $\chi(G)$ is inapproximable within
  $n^{1-\ve}$.
In our case $n \approx \sqrt{|V(G^+)|}$ holds, which yields a
 multiplicative error tending to infinity faster than
  $|V|^{1/2-\ve}$ in the approximation of $\cir(G^+)$
  as $n\to \infty$.
\qed

\bsk

We now turn to the positive result.
It requires a technical introduction before we state the theorem.

\msk

The treewidth of a graph $G$, denoted by $\tw(G)$, is equal to
 $\min \left(\omega(H)-1\right)$, where the minimum is taken over all
 \emph{chordal} graphs $H \supseteq G$.
From an algorithmic approach, treewidth equivalently is introduced
 via tree decompositions; we shall use a more specific kind of
 them as defined below.
For the fundamentals of the theory on treewidth, we refer to \cite{K94} and chapters 7 and 11 of \cite{8aut}.

Given any graph $G=(V,E)$, a \emph{nice tree decomposition}
 $\cT$ of $G$ consists of a rooted binary tree $T$ whose nodes will be denoted by
 $x_1,\dots,x_k$, together with non-empty subsets $V_1,\dots,V_k\subset V$
 where each node $x_i$ is associated with the corresponding $V_i$.

Two types of restrictions are put on the sets $V_i$.
One type with three conditions is related to $G$, namely
 \begin{itemize}
  \item[$(i)$] $V_1 \cup \cdots \cup V_k = V$;
  \item[$(ii)$] if $vw\in E$ then there is a node $x_i$ where
    $v,w\in V_i$;
  \item[$(iii)$] if $v\in V_{i'}$ and $v\in V_{i''}$ then also
    $v\in V_i$ holds for all $i$ such that $x_i$ is an internal
    node of the unique $x_{i'}$--$x_{i''}$ path in $T$.
 \end{itemize}
In order to have a clear distinction between the two structures, we use the term ``vertices'' in the graph $G$ and ``nodes'' in the host tree $T$ of its tree decomposition.

The other type of restrictions categorize the nodes $x_i$ in
 terms of their down-degree in $T$ and associated set $V_i$,
 as follows:
  \begin{itemize}
    \item a \emph{leaf node} $x_i$ has no children in $T$
    
    \item an \emph{introduce node} $x_i$ has one child $v_{i'}$
     in $T$, and its set $V_i$ is obtained from $V_{i'}$ by
     inserting just one vertex, i.e.\ $V_i = V_{i'}\cup\{v\}$
     for some $v\in V\setminus V_{i'}$;
    \item a \emph{forget node} $x_i$ has also one child $x_{i'}$
     in $T$, but its set $V_i$ is obtained from $V_{i'}$ by
     omitting just one vertex, i.e.\ $V_i = V_{i'}\setminus\{v\}$
     for some $v\in V_{i'}$;
    \item a \emph{join node} $x_i$ has two children
     $x_{i'},x_{i''}$ in $T$, and all their sets are the same,
     i.e.\ $V_i = V_{i'} = V_{i''}$.
  \end{itemize}
The \emph{width} of $\cT$ is
 $\displaystyle \max_{1\leq i\leq k} \left(|V_i|-1\right)$.
Theory proves that $\tw(G)\leq t$ holds if and only if
 $G$ admits a nice tree decomposition having width at most $t$,
 that is $|V_i|\leq t+1$ for all $i$.
It is also known that in this case the number of nodes in the
 host tree $T$ need not exceed $4|V|$, hence it can be ensured
 to be linear in the order of $G$.

 For later reference, we denote by $V_r$ the subset of $V$ associated with the root of $T$.

\bthm
For every positive integer $t$, the values of $\ali$, $\omi$, $\cir$, and $\thi$ can be determined in linear time on graphs of treewidth at most $t$.
\ethm

\pf
Let $\cG=\cG _t$ be the class of graphs $G$ with $\tw(G)\leq t$,
 for a fixed positive integer $t$.
Consider a generic input graph $G=(V,E)$ from $\cG$.
We take a nice tree decomposition $\cT=(T;V_1,\dots,V_k)$
 of width $t$ and $|V(T)|=O(|V|)$.
A dynamic programming algorithm will be applied along a
 postorder traversal of $T$.
For each node $x_i$ of $T$ a computational table $\tb_i$
 will be determined.

The indexing of rows in the tables will have two major parts.
The first part gives information about the 1-selection under
 consideration; this part is analogous in all the four problems
 $\ali,\omi,\cir,\thi$.
The second part is more problem-specific, as it will be detailed later.

For any $V_i$, the components of the first part of row indexing are:
 \begin{itemize}
  \item a partial (possibly empty) 1-selection $f$ of edges
   inside the subgraph $G[V_i]$ induced by $V_i$ in $G$;
  \item a partition $V_i^+\cup V_i^- = V_i$;
  \item the subset $V_i^+\subset V_i$ consists of those vertices
   for which   it is assumed that a 1-selection has already
   been made, either inside $V_i$ or with an edge whose other
   end is in the earlier (already forgotten) subgraph of $G$;
  \item the subset $V_i^-\subset V_i$ of vertices for which
   it is assumed that no 1-selection has been made yet.
 \end{itemize}
 We note that $V_r^+ = V_r$ can be assumed without loss of generality, but this is not the case at nodes different from the root.
 
Analogously to the concept of $D(G, f)$ proposed in the Introduction, it is convenient to represent the 1-selection $f$ inside $V_i$
 by a directed graph, where an arc $(v,w)$ means that the edge
 $vw$ of $G$ is assigned to $v$ by $f$.
This information can be handled in the tables of the four
 node types as follows.
 \begin{itemize}
  \item If $x_i$ is a leaf node, then every subset consisting of
   vertices non-isolated inside $V_i$
   has to be considered as a $V_i^+$, and for each $V_i^+$
   all possible 1-selections have to be taken in the indexing
   of rows of the table for $x_i$.
  \item If $x_i$ is an introduce node with child $x_{i'}$,
   and the new vertex in $V_i$ is $v$, then: the option $v\in V_i^-$
   has to be taken with all cases of $V_{i'}$; and if $v$ has
   at least one neighbor in $V_{i'}$, then also $v\in V_i^+$
   has to be considered with every possible 1-selection at $v$.
  Moreover, for the subset of vertices $v'\in V_{i'}^-$ that are
   adjacent to $v$, all combinations of edges for a 1-selection
   have to be taken; the corresponding vertices are then moved
   from $V_i^-$ to $V_i^+$.
  \item If $x_i$ is a forget node, then the ``forgotten'' vertex
   is just removed from $V_{i'}$, the status of the remaining
   vertices and 1-selection inside $V_i$, is unchanged.
  \item If $x_i$ is a join node, then it is necessary to check
   that the cases at $V_{i'}$ and $V_{i''}$ are compatible.
  This means not only that we have the same 1-selection $f$ and
   the same partition (i.e., $V_{i'}^-=V_{i''}^-=V_i^-$ and
   $V_{i'}^+=V_{i''}^+=V_i^+$) at the two children.
  If a vertex is in $V_i^+$, then its selected edge must be
   inside $V_i$.
 \end{itemize}
In the recursive computation of parameters, all the subgraphs
 obtained by deleting the 1-selections will be considered.

\msk

\underline{Computing $\omi$:}

\smallskip

A complete subgraph after the removal of a 1-selection is
 complete also in $G$, and its vertices appear together in
 at least one $V_i$.
So in each $V_i$ we register all possible complete subgraphs $K\subset G_f[V_i]$
 for every $f$ and every $(V_i^-,V_i^+)$,
  and compute a value $w(K)$.

At a leaf node, $w(K)$ is the number $|V(K)|$ of vertices in $K$.

At an introduce node with $V_i=V_{i'}\cup \{v\}$, $w(K)$ taken from the table of $V_{i'}$ is
 unchanged if $v\notin V(K)$, and otherwise it is
 $w(K) = \max \left( w(K-v), |V(K)| \right)$. 

At a forget node with $V_i=V_{i'}\setminus \{v\}$, $w(K)$ is
 redefined as the maximum of its former value at $x_{i'}$
 and that of $w(K\cup \{v\})$. 

At a join node, values $w(K)$ are available in the tables of its
 two children.
Then the updated $w(K)$ is the larger of the two.

Then $\omega_1(G)$ can be read from the root as follows: for every partial 1-selection $f$ in $G[V_r]$, one takes the maximum over all the $w$-values in rows corresponding to cliques $K\subseteq G[V_r]_f$, and then one takes the minimum over all $f$s.

\msk

\underline{Computing $\ali$:}

\smallskip

This algorithm is essentially the same as the one determining
 the independence number on graphs of bounded treewidth.
The difference is that the possible removals of 1-selections
 have to be taken into account, and the independent sets of
 those subgraphs are listed.

At a leaf node, all independent sets $S$ are listed, and the
 value is $w(S)=|S|$.

At an introduce node with $V_i=V_{i'}\cup \{v\}$, the value
 $w(S)$ remains unchanged if $v\notin S$, and it is computed as
 $w(S):=w(S-v)+1$ if $v\in S$.

At a forget node with $V_i=V_{i'}\setminus \{v\}$, the formula
 $w(S):=\max\left( w(S), w(S-v) \right)$ is applied.

At a join node, the value $w(S)$ is computed as the sum of the
 two $w(S)$ values at the children, minus $|S|$.

In the end, $\cir(G)$ is equal to the largest value of $w(S)=w(S,f)$ at the root of $T$, taken over all partial 1-selections $f$ in $G[V_r]$ and all independent sets $S\subseteq G[V_r]_f$.

\msk

\underline{Computing $\cir$:}

\smallskip

Since $\cir(G)\leq \chi(G)\leq \tw(G)+1$ holds for every graph
 $G$, we know that $\cir$ is bounded above by a constant.
Then a simple linear-time algorithm to test whether $\cir(G)\leq  k$
 holds is obtained by generating all proper $k$-colorings of $V_i$
 with respect to $f$, and checking which of them is compatible
 (also regarding the partition $(V_i^+,V_i^-)$)
 with at least one such coloring at each child node of $x_i$.

 In the end, $\cir(G)$ is equal to the smallest $k$ for which the algorithm above terminates with an admissible coloring at the root node.
 
\msk

\underline{Computing $\thi$:}

\smallskip

Here at each node $x_i$ for each $f$ and each $(V_i^+,V_i^-)$
 we need to generate all partitions $\cP$ of $V_i$ such that
 each partition class is a complete subgraph after the removal
 of the edges selected by $f$.
Moreover, it is necessary to distinguish between two possibilities
 for each complete subgraph selected as a class in $\cP$.
Namely, whether it is assumed to contain an already ``forgotten''
 vertex in the computation or it did not have any vertex outside
 $V_i$ previously.

At a leaf node, no class is associated with forgotten vertices,
 and the value $w(\cP)$ of $\cP$ is the number of its classes.

At an introduce node with $V_i=V_{i'}\cup \{v\}$, attaching $v$
 to a class of $V_{i'}$ is feasible only if no forgotten
  vertices are associated with that class.
If $v$ is attached to an existing class, then $w(\cP)$ remains
 the same as $w(\cP-v)$ in $V_{i'}$; otherwise, if $\{v\}$ is a
 new singleton class, then $w(\cP)=w(\cP-v)+1$.

At a forget node with $V_i=V_{i'}\setminus \{v\}$, the value of a
 partition does not change; but the status of the class from which
 $v$ has been removed will indicate from then on that it is
 associated with a forgotten vertex.

At a join node $x_i$, it is not allowed to keep a partition $\cP$
 if it has a class with associated forgotten vertices at both
 children of $x_i$.
(Apart from this condition, both children may associate
 forgotten vertices with any number of partition classes.)
If a partition $\cP$ is kept for $V_i$, then its value is the
 sum of values at the two children of $x_i$, minus the number of
 classes in $\cP$.

At the end, $\thi(G)$ is equal to the smallest value of $w(\cP)=w(\cP,f)$ at the root of $T$, taken over all partial 1-selections $f$ in $G[V_r]$, where $\cP$ is the trivial partition with $V_r^+ = V_r$.
\qed

\section{Concluding remarks}

This paper presents a systematic study of a new graphical invariant called the robust chromatic number, motivated by its applicability in extremal combinatorics. In addition, we introduce ``robust versions'' of several fundamental graph parameters, including the independence number, clique number, clique covering number, and chromatic index. Basic estimates and relationships to other parameters are established, and algorithmic aspects are also considered to some extent. While some of the new results parallel classical ones, others are distinct and unique.

One can naturally extend the robust version of any other graph invariant following the same approach used to obtain $\cir$ from $\chi$ or $\omi$ from $\omega$, etc. This opens up a promising new area for future research. Although we do not provide an explicit list of parameters here, we propose and encourage a systematic exploration of this aspect. In particular, any variant of graph coloring presents an interesting direction for further investigation.

Besides these very general suggestions, we list here some more
 definite problems that remain open in connection with the
 robust chromatic number.
The first question concerns a possible strengthening in
 part $(ii)$ of Theorem \ref{t:pla}.

\begin{prm}
Do there exist planar graphs with\/ $\cir(G)=3$, or is\/ $2$
 a universal upper bound?
\end{prm}

It is a well-known elementary fact that the chromatic number is
 additive with respect to the complete join operation.
This is not the case for $\cir$, as shown by many examples above.

\begin{prm}
$(i)$ \
Is there a transparent way to determine\/ $\cir(G \oplus H)$,
 at least if\/ $\cir(G)$ and\/ $\cir(H)$ are also given,
 possibly with optimal\/ $1$-selections\/ $f_G$ and\/ $f_H$\,?

$(ii)$ \
Is there a natural graph operation for which\/ $\cir$ is
 additive on vertex-disjoint graphs?

$(iii)$ \
Is there a natural analogue of the class of cographs (\/$=$
 the graphs not containing any induced\/ $P_4$ subgraph)
 for\/ $\cir$\,?
\end{prm}


There seems to be a lot to do in strengthening the estimates
 in part $(iv)$ of Theorem \ref{t:oper}
 for the union of $k$ graphs, where the currently available constructions are very limited.

\begin{prm}
$(i)$ \ 
Find matching lower and upper bounds on
 the robust chromatic number of the union of\/ $k$ graphs.

$(ii)$ \ 
Given two integers\/ $k,t\ge 2$, compare\/
$\cir(G_1\cup\cdots\cup G_k)$ with\/
 $\prod_{i=1}^k \cir(G_i)$ under the assumption that
 each\/ $G_i$ has\/ $\cir(G_i)\geq t$.
\end{prm}

The line graph operation seems to be of interest in its own right.

\begin{prm}
$(i)$ \
Describe further infinite classes of graphs whose members\/
 $G$ satisfy\/ the equality\/ $\cir'(G)=\cir(L(G))$.

$(ii)$ \
Does\/ $\cir'(G)\leq \cir(L(G))$ hold for every graph\/ $G$\,?
\end{prm}

So far very little is known about the complexity of determining the robust parameters of graphs.  $\iota$

\begin{prm}
$(i)$ \
Describe classes of well-structured graphs on which\/ $\cir$
 can be determined in polynomial time.

$(ii)$ \
Describe classes of well-structured graphs on which the
 computation of\/ $\cir$ is \nph.

$(iii)$ \
Study the analogous problems for the related graph invariants\/
 $\omi$, etc., introduced above.

$(iv)$ \
Describe conditions in terms of forbidden subgraphs and forbidden induced subgraphs,
 under which the computation of various robustness parameters becomes tractable.
\end{prm}

\begin{prm}
Study the properties of robust total coloring and its
 parameter\/ $\cir''$.
\end{prm}

\paragraph{Acknowledgements.}

This research was supported in part by the National Research, Development and Innovation Office --
NKFIH under the grants SNN 129364 and FK 132060.


\begin{thebibliography}{99}

\bibitem{8aut}
M. Cygan, F. V. Fomin, \L. Kowalik, D. Lokshtanov, D. Marx, M. Pilipczuk, M. Pilipczuk, S. Saurabh:
\textit{Parameterized Algorithms.}
Springer Cham, 2016.

\bibitem{H}
S. L. Hakimi:
On the degrees of the vertices of a directed graph.
\textit{J. Franklin Inst.} \textbf{279} (1965), 290--308.

\bibitem{K94}
T. Kloks: Treewidth, Computations, and Approximations.
\textit{Lecture Notes in Computer Science}, \textbf{842}, Springer, 1994.

\bibitem{L1} L. Lovász: Normal hypergraphs and the perfect graph conjecture. \textit{Discrete Math.} \textbf{2} (1972), 253--267.

\bibitem{L2} L. Lovász: A characterization of perfect graphs. \textit{J. Comb. Theory, Ser. B} \textbf{13} (1972), 95--98.

\bibitem{Padb} M. Padberg: Perfect zero-one matrices. \textit{Math. Program.} \textbf{6} (1974), 180--196.

\bibitem{PTV}
B. Patk\'os, Zs. Tuza, M. Vizer:
Extremal graph-theoretic questions for $q$-ary vectors.
Manuscript, 2022.

\bibitem{Sei} D. Seinsche: On a property of the class of $n$-colorable
graphs. \textit{J. Combin. Theory Ser. B} \textbf{16} (1974), 191--193.

\bibitem{STV}
M. Stiebitz, Zs. Tuza, M. Voigt:
Orientations of graphs with prescribed weighted out-degrees.
\textit{Graphs Combin.} \textbf{31} (2015), 265--280.

\bibitem{Z}
D. Zuckerman:
Linear degree extractors and the
inapproximability of Max Clique and Chromatic Number.
\textit{Theory Comput.} \textbf{3} (2007), 103--128.

\end{thebibliography}
\end{document}